\documentclass[12pt]{article}
\usepackage{anysize}
\usepackage{amsmath}
\usepackage{amsthm}
\usepackage{amsfonts}
\usepackage{epsfig}
\usepackage{url}
\DeclareMathOperator{\sd}{sd}
\DeclareMathOperator{\cd}{cd}
\newcommand{\KG}{\mathrm{KG}}
\DeclareMathOperator{\alt}{alt}
\DeclareMathOperator{\conn}{connectivity}
\DeclareMathOperator{\cone}{cone}
\DeclareMathOperator{\susp}{susp}
\newcommand{\ZZ}{{Z}}
\newcommand{\lra}{\longrightarrow}
\newcommand{\lmt}{\longmapsto}
\newcommand{\sse}{\subseteq}
\newcommand{\sm}{{\setminus}}
\renewcommand{\SS}{{\cal S}}
\newcommand{\TT}{{\cal T}}
\newcommand{\CC}{{\cal C}}
\newcommand{\pmz}{\{{+},{-},0\}}
\newcommand{\disj}{_{\rm{\kern-2.7pt disj}}}
\newcommand{\stab}{_{\rm{\kern-2.7pt stab}}}
\newcommand{\tstab}{_{\rm{\kern-2.7pt {\it t\/\textrm-}stab}}}
\newcommand{\rstab}{_{\rm{\kern-2.7pt {\it r\/\textrm-}stab}}}
\newcommand{\chizero}{K} 
\newcommand{\sups}{^{\mbox{\itshape\bfseries\scriptsize s}}}
\newcommand{\s}{\mbox{\itshape\bfseries s}}
\newcommand{\ds}{{\s}}
\newcommand{\subds}{_{\mbox{\itshape\bfseries\scriptsize s}}}
\newcommand{\ssum}{{\overline n}}
\def\K{{\hbox{\sf K}}}
\def\L{{\hbox{\sf L}}}
\def\lsoft{{l\kern-0.035cm\char39\kern-0.03truecm}}
\newcommand\Topologically{\begin{small}\begin{quote}{\itshape{}Topologically}}
\newcommand\TopologicallyEnd{\end{quote}\end{small}}
\theoremstyle{plain}
\newtheorem{thm}{Theorem}[section]
\newtheorem{cor}[thm]{Corollary}

\newtheorem{lemma}[thm]{Lemma}

\theoremstyle{definition}
\newtheorem{defn}[thm]{Definition}
\newtheorem{example}[thm]{Example}
\theoremstyle{remark}

\title{Generalized Kneser Coloring Theorems\\
with Combinatorial Proofs}
\author
{{\Large G\"unter M.~Ziegler%
\footnote{%
Supported by Deutsche Forschungs-Gemeinschaft (DFG), %
and by the Miller Institute at UC Berkeley}}\\
Department of Mathematics, MA 6-2\\
TU Berlin, D-10623 Berlin, Germany\\
\url{ziegler@math.tu-berlin.de}}
\begin{document}
\maketitle

\begin{abstract}\noindent
The Kneser conjecture (1955) was proved by Lov\'asz (1978)
using the Borsuk-Ulam theorem; all subsequent proofs,
extensions and generalizations
also relied on Algebraic Topology results, 
namely the Borsuk-Ulam theorem and its extensions.
Only in 2000, Matou\v sek provided the first combinatorial
proof of the Kneser conjecture.

Here we provide a hypergraph coloring theorem,
with a combinatorial proof, 
which has as special cases the Kneser conjecture
as well as its extensions and generalization by
(hyper)graph coloring theorems of 
Do\lsoft{}nikov, Alon-Frankl-Lov\'asz, Sarkaria, and Kriz.
We also give a combinatorial proof of Schrijver's theorem.
\end{abstract}

\section{Introduction}

Kneser's conjecture \cite{Kne} stated that every coloring
of the graph $\KG^2{[n]\choose k}$, which 
has vertex set $[n]\choose k$, and in which 
two $k$-subsets of an $n$-set are connected by an edge if they are disjoint,
needs $n-2k+2$ colors for a proper vertex coloring, if $n\ge2k\ge4$.

Kneser's conjecture was first proved by 
Lov\'asz \cite{Lo}, in one of the first, and most 
spectacular, applications of an Algebraic Topology result
(the Borsuk-Ulam theorem) to a combinatorial problem.
An alternative proof was later given by B\'ar\'any \cite{Ba}, 
extensions by Schrijver \cite{Schrijver}, 
Alon, Frankl \& Lov\'asz \cite{AFL}, Do{\lsoft}nikov \cite{Dol1},
Sarkaria \cite{Sa1}, and finally by 
Kriz \cite{Kri1,Kri2}, whose result implies the
theorems by Lov\'asz, Do{\lsoft}nikov and Alon-Frankl-Lov\'asz.
All of these were proved using
Algebraic Topology results, namely the Borsuk-Ulam theorem
and its extensions. 
This established them as a part of the
``classical core'' of Topological Combinatorics \cite{Bj} \cite{BMZ}.

A curious aspect of this is that the basic
topological result used in this context, the 
Borsuk-Ulam theorem, has a variety of ``combinatorial proofs,''
that is, reductions via simplicial approximation to combinatorial 
results such as Tucker's lemma \cite{Tuck},
the Ky Fan lemma \cite{Fan}, etc. 

In 2000, Ji\v{r}\'{\i} Matou\v sek provided two breakthroughs 
for this situation. 
\emph{First}, in \cite{Mat1} he provided
a combinatorial bypass of the Borsuk-Ulam theorem, and thus
a combinatorial proof of Kneser's conjecture. This used
only an entirely combinatorial special case of 
the Tucker lemma.
\emph{Secondly}, in \cite{Mat2}, Matou\v{s}ek gave a simple and elegant
derivation of Kriz' theorem from Dold's theorem.

Here, we shall demonstrate the power and extend the scope of Matou\v sek's 
approach, by establishing
\begin{itemize}\itemsep=-4pt
\item a simple combinatorial proof of Do{\lsoft}nikov's theorem
(an extension of Matou\v sek's proof \cite{Mat1} of the Kneser conjecture),
\item a new, fairly general hypergraph coloring theorem 
(which has the theorems by Lov\'asz \cite{Lo}, 
Alon, Frankl \& Lov\'asz \cite{AFL}, Do{\lsoft}nikov \cite{Dol1},
Sarkaria \cite{Sa1}, and Kriz \cite{Kri1,Kri2} as special cases),
together with a combinatorial proof, and
\item a combinatorial proof of Schrijver's theorem
(via cyclic oriented matroids).
\end{itemize}
The proofs that we give are combinatorial (``elementary'') in the sense that 
they do not rely on topological concepts (such as continuous maps, simplicial 
approximations, homology) or results. We do, however, 
phrase some proofs in terms of
simplicial complexes and maps, 
chain complexes, etc.: 
their use can be eliminated at the price of making the phrasing of
the proofs more cumbersome. 
More importantly, topological knowledge, interpretation and 
inspiration is ``behind'' the proofs of this paper: it is neither
desirable nor practical to eliminate this background intuition.

\section{Some Tools and Notation}

We use $[n]$ to denote the set $\{1,\ldots,n\}$ of~$n$  
integers, and $2^N$ for the set of all subsets of a finite set~$N$.
A \emph{hypergraph} is a set family $\SS\sse2^N$; the set $N$ is its 
\emph{ground set}.
The hypergraphs that appear in the following will have \emph{no loops},
that is, all their edges $S\in\SS$ have cardinality at least~$2$.
A hypergraph is \emph{$r$-uniform} if all of its \emph{edges} $S\in\SS$ 
have the same cardinality~$r$.
The \emph{restriction} of a hypergraph $\SS\sse2^N$ to a subset
$M\sse N$ of its ground set is the hypergraph 
$\SS|_M=\{S\in\SS:S\sse M\}\sse 2^M$.

In the following $\s=(s_1,...,s_n)$ will be a vector
of \emph{multiplicities}, with $1\le s_i<r$. The sum of multiplicities
will be denoted by $\ssum:=s_1+\ldots+s_n$. Usually 
$\s$ will be \emph{constant}, that is, $s_1=\ldots=s_n$, with $\ssum=ns$.
An $r$-tuple of subsets $R_1,\ldots,R_r$
is \emph{$\s$-disjoint} if each $i\in[n]$ is contained in
at most $s_i$ of the sets $R_j$, that is,
if $|\{j\in[r]:i\in R_j\}|\le s_i$ for all~$i$.
For example, an $(1,...,1)$-disjoint family of sets is simply disjoint;
in this case $\ssum=n$.
An $(s,\ldots,s)$-disjoint family is what Sarkaria \cite{Sa1}
calls ``$(s+1)$-wise disjoint.''

Interesting examples of hypergraphs that appear in the following include
$[n]\choose k$, the collection of all subsets of~$[n]$ of cardinality~$k$,
and ${[n]\choose k}\stab$, the collection of all \emph{stable $k$-subsets},
that is, all subsets that do not contain any two adjacent elements modulo~$n$. 
More generally, a subset $S\sse[n]$ is \emph{$t$-stable} if
any two of its elements are at least ``distance~$t$ apart'' on
the $n$-cycle, that is, if $t\le |i-j|\le n-t$ for distinct $i,j\in S$.
Thus every set $R\sse [n]$ is $1$-stable, while
stable is the same as $2$-stable.

A \emph{coloring} of a hypergraph $\SS\sse2^N$ with $m$ \emph{colors}
is a function $c:N\lra[m]$ that assigns colors to the ground set
so that no edge $S\in\SS$ is monochromatic, that is, every edge
contains two elements $i,j\in S$ with $c(i)\neq c(j)$.
Equivalently, no $c^{-1}(i)$ contains a set $S\in\SS$.
The \emph{chromatic number} $\chi(\SS)$ of a hypergraph is
the smallest number $m$ such that a coloring of~$\SS$ with $m$ colors
exists. (Clearly $\chi(\SS)\le|N|$ if $\SS$ has no loops.) 

The \emph{$r$-colorability defect} of~$\SS\sse2^{N}$ is the number of
elements one has to remove from the ground set of~$\SS$ so that
the remaining hypergraph can be colored with $r$ colors, that is,
the smallest cardinality of a subset $N\sm M$ of~$N$ such that
$\chi(\SS|_M)\le r$.
(This is Matou\v{s}ek's terminology for a concept 
introduced by Do\lsoft{}nikov and Kriz.)
More generally, the \emph{$\s$-disjoint $r$-colorability defect} 
of~$\SS\sse2^N$ is
\\$
\cd^r\subds\SS:= 
\ssum-\max\big\{\sum\limits_{j=1}^r|R_j|:
R_1,\ldots,R_r\sse N \mbox{ $\s$-disjoint\,},\ 
S\not\sse R_j\mbox{ for }S\in\SS\big\},
$\\
that is, the number of elements that have to be removed from
the multiset $N\sups$ so that the remaining multiset can be covered
by an $\s$-disjoint $r$-family of sets (``color classes'') 
such that none of the sets
contains a set from~$\SS$ (``there is no monochromatic $\SS$-set'').

The \emph{$r$-th Kneser hypergraph} $\KG^r\SS$ of~$\SS$ is the $r$-uniform
hypergraph with ground set~$\SS$ whose edges are formed by 
the $r$-sets of disjoint edges of~$\SS$.
For example, if $\SS\sse{N\choose2}$ is a graph, then
the edges of $\KG^r\SS$ correspond to the matchings
in~$\SS$ of size~$r$.
More generally, the \emph{$r$-th $\s$-disjoint Kneser hypergraph}
of~$\SS$, denoted $\KG^r\subds\SS$,
is the $r$-uniform hypergraph with ground set~$\SS$ whose edges are 
the $\s$-disjoint $r$-sets of edges in~$\SS$.
\medskip

We write $\pmz^N$ for the set of all \emph{signed subsets} of~$N$,
that is, the family of all pairs $(X^+,X^-)$ of disjoint subsets
of~$N$. For $N=[n]$, such subsets can alternatively be encoded by 
\emph{sign vectors} $X\in\pmz^n$, where $X_i={+}$ denotes that $i\in X^+$,
while $X_j={-}$ means that $j\in X^{-}$. 
The \emph{positive part} of~$X$ is $X^+:=\{i\in[n]:X_i={+}\}$, 
and analogously for the \emph{negative part}~$X^-$.
In the following, we shall
switch freely between the different notations for signed sets.
For sign vectors, we use the usual partial order from
oriented matroid theory \cite{BLSWZ}, which is defined
componentwise with $0\le{+}$ and $0\le{-}$. Thus $X\le Y$, that is
$(X^+,X^-)\le(Y^+,Y^-)$, holds if and only if $X^+\sse Y^+$ and $X^-\sse Y^-$.
We write $\alt(X)$ for the
length of the longest alternating subsequence of non-zero signs in~$X$.
For example, $\alt(0{+}{-}{+}{+}00{+}0)=3$,
while $\alt(0{+}{-}{+}{+}{-}0{+}0)=5$.

$\ZZ_p$ denotes a cyclic group of order~$p$.
We interpret it as the group of $p$-th complex roots of unity, 
$\ZZ_p=\{1=\omega^p,\omega,\ldots,\omega^{p-1}\}$,
and so its elements are called \emph{signs}.
This will below be used for more general ``signed sets,'' where
elements get signs from~$Z_p$.

The Borsuk-Ulam theorem 
asserts that there is no $\ZZ_2$-equivariant (continuous) map from $S^d$ to
$S^{d-1}$. Dold's theorem \cite{Dold} is a transformation group
extension of this: for every $\ZZ_p$-equivariant map
$f:X\lra Y$ between free $\ZZ_p$-spaces (compact CW complexes, say)
the dimension of~$Y$ is larger than the connectivity of~$X$.
\medskip

The following sketch of proof for Dold's theorem (following
\cite{Dold}) is a ``blue\-print'' for the combinatorial
proofs in the following. Assume that $f: X\lra Y$ is a simplicial
or cellular map. If $\conn(X)\ge \dim(Y)$, then one can construct
an equivariant map back, $g: Y\lra X$, whose image is contained in a cone.
Since $g f$ is an equivariant map, its Lefschetz number
$\Lambda(g f)$, which \emph{counts} the simplices/\allowbreak cells
that are mapped to themselves (with signs according to dimension
and orientation reversal), is divisible by~$p$.
On the other hand, one can restrict $g f$ to a cone in~$X$, 
and from this derive that $\Lambda(g f)=1$, a contradiction.

For combinatorial proofs, the hard work is usually in the explicit
construction of the map back, $g: Y\lra X$, without use of, or
reference to, connectivity information. The Lefschetz number of a
chain map on a cone is dealt with by Lemma~\ref{l:cone}.

\section{Colorings and Colorability Defects}

For $\s=(1,...,1)$, the following coloring of the Kneser hypergraphs 
is due to Kneser~\cite{Kne} in the case $r=2$ 
and to Erd\H os~\cite{Er} in the general case. 
It corrects the coloring given in \cite[(3.3)]{Sa1}.

\begin{lemma}\label{l:col} 
For $r\ge2$, $k\ge2$, constant $\s=(s,\ldots,s)$ 
with $1\le s<r$, and $sn\ge kr$,
\[
\chi\big(\KG^r\subds\textstyle{[n]\choose k}\big)\ \ \le\ \ 
1+\big\lceil\tfrac1{\left\lfloor\tfrac{r-1}s\right\rfloor}\tfrac{ns-rk+1}s\big\rceil.
\]
\end{lemma}

\begin{proof}
Set $P:=\lfloor\frac{r-1}s\rfloor$ and 
$M:=\lceil\frac1P\frac{ns-rk+1}s\rceil$. 
With this an explicit coloring is given by
\[
S\ \ \lmt\ \ \min\big\{\big\lceil\tfrac1P\min(S)\big\rceil
\,,\,M+1\big\}.
\]
This rule assigns to each $k$-set $S$ an integer between $1$ and $M+1$.

If $\{S_1,\ldots,S_r\}$ is an $\s$-disjoint $r$-family,
then every minimal element $\min(S_j)$ can appear at most $s$ times
in the family; thus if $\lceil\frac1P\min(S_j)\rceil\le M$,
then this value 
is assigned to at most $Ps$ sets in an $\s$-disjoint $r$-family, where
$Ps=\lfloor\frac{r-1}s\rfloor s\le r-1$. 

On the other hand, if all the $k$-sets $S_j$ get the color $M+1$,
then they are contained in the set $\{PM+1,\ldots,n\}$, of cardinality 
$n-PM=n-P\lceil\frac1P\frac{ns-rk+1}s\rceil \le n-P\frac1P\frac{ns-rk+1}s=
       \frac{rk-1}s<\frac{rk}s$.
But the pigeonhole principle demands that an $s$-disjoint $r$-family of
$k$-sets uses at least $\frac{rk}s$ elements.
\end{proof}

The colorings of Lemma~\ref{l:col} will be shown to be
optimal whenever $s$ divides $r-1$:
see Section~\ref{sec:badcolor}, where we also analyze
a case where the coloring is far from optimal.

Lemma~\ref{l:col} also provides colorings for the induced sub-hypergraphs
$\KG^r{[n]\choose k}\tstab\sse \KG^r{[n]\choose k}$, for $t\ge1$, so
$\chi\big(\KG^r{[n]\choose k}\tstab\big)\le\lceil\tfrac{n-(k-1)r}{r-1}\rceil$.
This coloring is still optimal for $r=t=2$: this is Schrijver's theorem
(see Section~\ref{sec:schrijver}).
For $r>2$, see Section~\ref{s:conj:big}.
\medskip

The theorems by Do\lsoft nikov and Kriz and our Theorem~\ref{t:main}
give lower bounds for chromatic numbers of hypergraphs in terms
of $r$-colorability defects. These lower bounds are useful 
only since they are easy to evaluate, e.~g.\ as follows.

\begin{lemma}\label{l:cdefects} 
Let $r\ge2$, $n\ge k\ge 2$, $t\ge1$, 
and constant $\s=(s,...,s)$ with $1\le s<r$.
If~$n\ge tk$ (otherwise ${[n]\choose k}\tstab=\emptyset$), then 
\begin{eqnarray*}
\cd^r\subds\textstyle{[n]\choose k}\tstab 
 &=& \max\{ ns-tr(k-1),0\}.\\
\noalign{\noindent In particular $(t=1)$,}
\cd^r\subds\textstyle{[n]\choose k}\phantom{\tstab}
 &=& \max\{ ns-\,r\,(k-1),0\},\\
\noalign{\noindent and $(t=2)$}
\cd^r\subds\textstyle{[n]\choose k}\stab\ \ 
 &=& \max\{ ns-2r(k-1),0\}.
\end{eqnarray*}
\end{lemma}

\begin{proof}
No set $R\sse[n]$ of $t(k-1)$ contiguous elements mod~$n$
contains a $t$-stable $k$-set. (Note $t(k-1)<n$.)
Furthermore, there is an $\s$-disjoint packing
of (at most $ns$ elements from) $r$~such contiguous subsets into $[n]$:
Such a packing can be written down as
\\$
R_j\ \ :=\ \ 
\{(j-1)t(k-1)+1\mbox{\rm\ mod\,}n,\ 
  (j-1)t(k-1)+2\mbox{\rm\ mod\,}n,\ \ldots\ ,
      jt(k-1)  \mbox{\rm\ mod\,}n\}.
$\\
This proves that 
$\cd^r\subds{[n]\choose k}\tstab\le \max\{ ns-tr(k-1),0\}$.

To prove that 
$\cd^r\subds{[n]\choose k}\tstab\ge \max\{ ns-tr(k-1),0\}$,
it suffices to verify that every set $R\sse[n]$ of cardinality
$t(k-1)+1$ contains a $t$-stable $k$-subset.
Take $R$, and let $S\sse[N]\sm R$ be an arbitrary set of
size~$t-1$; this exists since $n\ge tk$. Now
$R\cup S$ has cardinality $tk$, and we can partition it
into $t$ disjoint $t$-stable $k$-subsets, by taking
``every $t$-th element'' to go into the same $k$-subset.
At least one of these $t$ $k$-subsets contains no element from~$S$,
since $|S|<t$.
\end{proof}

\section{Do{\lsoft}nikov's Theorem}

Tucker's lemma \cite{Tuck} says that if we take a suitable triangulation
of an $n$-ball, and label its vertices by labels in 
$\{\pm1,\ldots,\pm n\}$ in a way that is antipodal on 
the boundary, then there is a ``complementary edge''
whose endpoints receive opposite labels~$\pm i$.
Matou\v sek's proof \cite{Mat1} of the Kneser conjecture
relies on  the following combinatorial lemma,
which corresponds to Tucker's lemma applied to 
(the boundary of)
the barycentric subdivision of the $n$-cube, $\sd([-1,+1]^n)$,
whose vertex set can be identified with~$\pmz^n$.

\begin{lemma}[Octahedral Tucker lemma]\label{l:octa}
If $\lambda{:}\ \pmz^n\sm\{0\}^n{\longrightarrow}\{\pm 1,\ldots,\pm (n-1)\}$ 
satisfies $\lambda(-X)=-\lambda(X)$ for all $X$,
then there are signed sets $(A^+,A^-)$ and $(B^+,B^-)$
such that $\lambda(A^+, A^-) = - \lambda(B^+,B^-)$,
with $A^+\sse B^+$ and $A^-\sse B^-$.
\end{lemma}

\noindent
This lemma has simple combinatorial proofs, 
e.~g.\ by the method of Freund \& Todd \cite{FreundTodd}; see \cite{Mat1}.
(For further combinatorial Tucker lemmas, see Aigner~\cite{Ai}.)

\begin{thm}[Do{\lsoft}nikov \cite{Dol1}] \label{t:dol}
For every hypergraph $\SS\sse 2^{[n]}$,
the $2$-colorability defect is a lower bound for
the chromatic number,
\[
\chi(\KG^2\SS)\ \ \ge\ \ \cd^2\SS.
\]
\end{thm}

\begin{proof}[\bfseries Combinatorial Proof. ]
Let $c:\SS \longrightarrow [m]$
be a proper $m$-coloring, and assume that
$\cd^2\SS>m$, that is, 
if any subset of~$[n]$ of size at least $n-m$ is colored by two colors,
then it contains a monochromatic subset from $\SS$.
Fix an arbitrary linear ordering $\prec$ on the subsets of~$[n]$. 
Then define a map
$\lambda: \pmz^n\sm\{0\}^n \lra \{\pm 1,\ldots,\pm (n-1)\}$, as follows:
\begin{enumerate}\itemsep=-1pt
\item
If $|A^+| + |A^-| \ge n-m$, then define $\lambda(A^+ , A^-)$ as
$ \pm c(S)$, where $S$ is the smallest set (according to ``$\prec$'')
from $\cal S$ that is contained either in $A^+$, or in $A^-$.
Take the sign to indicate which of~$A^+$ or $A^-$ you took $S$ from.
Thus we obtain a value $\lambda(A^+ , A^-)$ in the set
$\{\pm 1, \pm2, \ldots, \pm m\}$.
\item
If $|A^+| +|A^-| \le n-m-1$, then define $\lambda(A^+ , A^-)$ as
$ \pm (m+ |A^+ |+| A^-|)$, where the sign indicates 
which of~$A^+$ or $A^-$ is nonempty, and if they both are,
then it indicates which is smaller (according to ``$\prec$'').
Thus we obtain a value $\lambda(A^+ , A^-)$ in the set
$\{\pm (m+1), \ldots, \pm (n-1)\}$.
\end{enumerate}
This map $\lambda$ is antipodal.
Thus by the Octahedral Tucker lemma~\ref{l:octa}, 
there are signed sets $(A^+, A^-)$, $(B^+, B^-)$
with $\lambda(A^+, A^-) = - \lambda(B^+, B^-) =\pm i$,
where $A^+\sse B^+$, $A^-\sse B^-$, not equality in both
cases, and so $|A^+ \cup A^-| < |B^+ \cup B^-|$.
This is possible only if both signed sets are labeled
according to the first case.
But then (assume without loss of generality that above we have ``$+i$'')
there are sets  $S,T\in\SS$ with $c(S)=c(T)=i$
and $S\sse A^+\sse B^+$, $T\sse B^-$, where $B^+$ and $B^-$ 
are disjoint: so also $S$ and $T$ are disjoint,
but they get the same color from~$c$, contradiction.
\end{proof}


\section{A Hypergraph Coloring Theorem}


\begin{thm}\label{t:main}
For every hypergraph $\SS\sse 2^{[n]}$, for $r\ge2$,
and for multiplicities $\s=(s_1,...,s_n)$ with $1\le s_i<r$,
the $\s$-disjoint $r$-colora\-bility defect yields a lower bound for
the chromatic number of the associated $r$-th $\s$-disjoint
Kneser hypergraph,
\[
\chi(\KG^r\subds\SS)\ \ \ge\ \ 
\left\lceil\tfrac1{r-1}\, \cd^r\subds\SS \right\rceil.
\]
\end{thm}

\noindent
This theorem, in combination with Lemma \ref{l:cdefects}, has many well-known
special cases, for constant $\s$:
\\
$\SS={[n]\choose k}$, $r=2$, $\s=(1,...,1)$: \ 
Lov\'asz \cite{Lo} (the Kneser conjecture),\\
$r=2$, $\s=(1,...,1)$: \ 
Do\lsoft nikov \cite{Dol1} (Theorem~\ref{t:dol}),\\
$\SS={[n]\choose k}$, $\s=(1,...,1)$: \ 
Alon, Frankl \& Lov\'asz \cite{AFL},\\
$\s=(1,...,1)$: \ 
Kriz \cite{Kri1,Kri2}, and\\
$\SS={[n]\choose k}$: \ 
Sarkaria \cite{Sa1}.
\\
The generalization
to non-constant $\s$ is not done for it's own interest,
but since it is needed for the first part of our proof,
where we show that one may assume that $\ssum-1$ is divisible by $r-1$. 
Under this assumption, and if $p:=r\ge2$ is a prime,
the second part of the proof derives the theorem
from the ``$\ZZ_p$-Tucker lemma'' \ref{l:Zp-Tucker} 
(this is where ``the topology is hidden''). 
The third part reduces the general case of the theorem to the prime case. 
Finally, in Section~\ref{Sec:Zp-Tucker},
the $\ZZ_p$-Tucker lemma is proved combinatorially.

\begin{proof}[\bfseries Reduction of Theorem \ref{t:main}
to the case when \boldmath$r-1$ divides $\ssum-1$.]
For this, we watch what happens if we increase the ground set,
by extending $[n]$ to~$[n{+}1]$, with \mbox{$s_{n+1}:=1$},
where $\SS\sse 2^{[n]}\sse2^{[n+1]}$ is not changed. 
Since $\SS$ is not changed,
the Kneser hypergraph $\KG^r\subds\SS$ and its chromatic 
number don't change, either.
On the other hand, with this operation $\ssum-1$ increases by~$1$, and 
\\$
\max\big\{\sum\limits_{j=1}^r |R_j|:\mbox{$\s$-disjoint $r$-family
$\{R_1,\ldots,R_r\}$ with no $S\not\sse R_j$}\big\}
$\\
also increases by~$1$, since we may extend exactly one
of the $R_j$s by an extra element $n+1$.
Thus, in summary, extending the ground set with $s_{n+1}=1$ changes
neither the chromatic number of the Kneser hypergraph,
nor the colorability defect, so validity of the theorem is
unchanged. By applying this operation, which increases
$\ssum-1$ by $1$, at most $r-2$ times, we
get the required divisibility.
\end{proof}

We write $\sigma^{n-1}$ for the $(n-1)$-dimensional
simplex with vertex set~$[n]$: this corresponds to the 
set system of faces~$2^{[n]}$.
Further, $\sigma^{n-1}_{k-1}$ denotes the
$(k-1)$-dimensional skeleton of this simplex, which
corresponds to the set system~$[n]\choose\le k$.

\begin{defn}[\boldmath$\ds$-disjoint $p$-fold joins]
If $\K$ is any simplicial complex on the 
ground set~$[n]$, then $\K^{*p}$ is the join of $p$ disjoint
copies of~$\K$, which is a simplicial 
complex on the ground set $\ZZ_p\times[n]$; this complex
has a natural $\ZZ_p$-action.

Similarly, the \emph{$\ds$-disjoint $p$-fold join} $\K^{*p}\subds$
is the complex of all subsets $A$ of the
ground set~$\ZZ_p\times[n]$, such that the elements with the same
``sign'' $\omega^k\in\ZZ_p$ correspond to a simplex in~$\K$, and such
that every element $i\in[n]$ appears in~$A$ with at most
$s_i$ different signs $\omega^k$.
This complex again has a natural $\ZZ_p$-action. 
If $p$ is prime, and if $s_i<p$ for all~$i$, then the
$\ZZ_p$-action on $\K^{*p}\subds$ is free.

We identify the ground set of~$\K^{*p}\subds$ with the
index set of an $n\times p$ matrix --- that is, with an
$n\times p$ \emph{chessboard} in the terminology of~\cite{BLVZ}.
Thus the faces of $\K^{*p}\subds$ may be viewed as
$0/1$-matrices of size $n\times p$, where
\begin{itemize}\itemsep=-3pt
\item in each column, the rows that contain a $1$ correspond to a face of~$\K$,
\item the $i$-th row contains at most $s_i$ ones, and
\item the $\ZZ_p$-action cyclically permutes the columns of the matrix.
\end{itemize}
The inclusion relation on faces of $\K^{*p}\subds$ translates
into the componentwise $\le$-partial order on $0/1$-matrices.
We write these matrices column-wise as $A=(A_1,\ldots,A_p)$,
where each~$A_j$ is the characteristic vector of a face of~$\K$.
\end{defn}

\begin{lemma}[\boldmath$\ZZ_p$-Tucker lemma] \label{l:Zp-Tucker}
Let $p\ge2$ be a prime, $n\ge1$, $\s=(s_1,...,s_n)$ with $1\le s_i<p$,
and let 
\begin{eqnarray*}
\lambda: \qquad 
(\sigma^{n-1})^{*p}\subds\sm \{\emptyset\}&\lra&\quad \ZZ_p\times[m] \\
         A\ =\ (A_1,\ldots,A_p)&\lmt&(\lambda_1(A),\lambda_2(A))=\lambda(A)
\end{eqnarray*}
be a $\ZZ_p$-equivariant map from non-zero faces/matrices in
$(\sigma^{n-1})^{*p}\subds$ to signed integers.

If $m\le \lfloor\frac{\ssum-1}{p-1}\rfloor$, then there is a 
chain of faces/matrices
\[
A^{(1)}\subset A^{(2)}\subset\ldots\subset A^{(p)}
\]
with $\lambda(A^{(i)}) = \omega^{\pi(i)}\lambda_2 (A^{(p)})$
for some permutation $\pi\in\Pi_p$, that is, such that
the~$A^{(i)}$
get assigned to the same absolute value $\lambda_2(A^{(i)})$,
but with $p$ distinct signs $\lambda_1(A^{(i)})\in\ZZ_p$.
\end{lemma}

\Topologically, 
this $\ZZ_p$-Tucker lemma can be derived from Dold's theorem:
If the conclusion does not hold, then $\lambda$ 
defines a $\ZZ_p$-equivariant simplicial map
\[
\lambda:\qquad
\sd\ (\sigma^{n-1})^{*p}\subds \ \ \lra\ \ 
     (\sigma^{m-1})^{*p}_{(p-1,...,p-1)}
\ =\  (\sigma^{p-1}_{p-2})^{*m}
\]
from the barycentric subdivision of the complex 
of all $0/1$-matrices of size $n\times p$ with at most $s_i$ ones 
in the $i$-th row, 
to the complex of all $0/1$-matrices of size $m\times p$ with at
most $p-1$ ones per row; this space can be written in two
different ways, depending on whether it is read
``column-wise'' as a deleted join, or ``row-wise''
as a proper join. 
On both spaces, the group $\ZZ_p$ acts by cyclic permutation of
the columns. If $p$ is prime and $s_i<p$, then the $\ZZ_p$-actions are free.

A maximal face of $(\sigma^{n-1})^{*p}\subds$
has exactly $s_i$ ones in the $i$-th row, so
the complex has dimension $\ssum-1$. We 
write the complex as $\sigma^{p-1}_{s_1-1} * \ldots * \sigma^{p-1}_{s_n-1} $
to conclude from the connectivity lemma for joins that
the connectivity of this complex 
(and of its barycentric subdivision) is $\ssum-2$.

The complex $(\sigma^{p-1}_{p-2})^{*m}$ is pure of dimension
$m(p-1)-1$; we don't even need that it is a simplicial sphere.
Its $\ZZ_p$-action is free since $p$ is a prime.
Thus we have a contradiction to Dold's theorem if $\ssum-2\ge m(p-1)-1$,
that is, if $m\le\lfloor\frac{\ssum-1}{p-1}\rfloor$. 
\TopologicallyEnd

\begin{proof}[\bfseries Proof of Theorem \ref{t:main} for prime~\boldmath$p$
and integral $\frac{\ssum-1}{r-1}$.]
Let $\chi(\KG^p\subds\SS)=K$, and let
$c:\SS \longrightarrow [K]$ be a coloring
such that no $p$ $\s$-disjoint sets from~$\SS$ get the same color.
At the same time we assume that $\cd^p\subds\SS>(p-1)K$, that is, if 
$\ssum-(p-1)K$ elements of the multiset $[n]\sups$
are colored by $p$ colors (which we take from~$\ZZ_p$), then 
some set from~$\SS$ is monochromatic. We define a labeling
\[
\lambda: \qquad(\sigma^{n-1})^{*p}\subds\sm\{\emptyset\}\ \ \longrightarrow\ \ 
 \ZZ_p\times \big[\textstyle{\lceil\frac{\ssum-1}{p-1}}\rceil\big],
\]
using an arbitrary linear ordering $\prec$ on the subsets of~$[n]$, as follows:
\begin{enumerate}\itemsep=-1pt
\item
If $|A_1|+\ldots+|A_p|\ge\ssum-(p-1)K$, then define $\lambda_2(A)$ as~$c(S)$, 
where $S$ is the smallest set (according to~``$\prec$'')
from $\cal S$ that is contained in one of the $A_i$'s;
take the sign $\lambda_1(A):=\omega^i\in \ZZ_p$ 
to indicate which $A_i$ you took $S$ from.
Thus we obtain a value $\lambda(A)=(\lambda_1(A),\lambda_2(A))$ 
in the set $\ZZ_p\times [K]$.
\item
If $|A_1|+ \ldots+ |A_p| \le \ssum-(p-1)K-1$, then define 
$\lambda_2(A):=K+\lceil\frac{|A_1|+ \ldots+ |A_p|}{p-1}\rceil$,
where the sign $\lambda_1(A):= \omega^i$ indicates which of the 
nonempty sets $A_i$ is the smallest one according to~``$\prec$''.
In this case we obtain a value $\lambda_2(A)$ in the set
$\{K{+}1, \ldots, \lceil\frac{\ssum-1}{p-1}\rceil\}$.
\end{enumerate}
This labeling is $\ZZ_p$-equivariant.
By assumption $\frac{\ssum-1}{p-1}$ is an integer,
that is, $\lceil \frac{\ssum-1}{p-1}\rceil=\lfloor\frac{\ssum-1}{p-1}\rfloor$.
Thus we can apply the $\ZZ_p$-Tucker lemma \ref{l:Zp-Tucker}: 
there is a chain of~$p$ $0/1$-matrices
$A^{(1)}\le A^{(2)}\le\ldots\le A^{(p)}$
such that $\lambda(A^{(i)})=(\omega^{\pi(i)},\lambda_2 (A^{(p)}))$,
for some permutation $\pi\in\Pi_p$.

Since at most $p-1$ of the matrices $A^{(i)}$ can have the same 
$\lceil\frac{|A^{(i)}_1|+\ldots+ |A^{(i)}_p|}{p-1}\rceil$,
and thus the same ``color'' $\lambda_2(A^{(i)})=k_0>K$ according to the
second case, the chain consists of~$p$ matrices
that fall into the first case in the definition of~$\lambda$.
Thus there are sets $S_i\in\SS$ that satisfy 
$S_i\sse A^{(i)}_{\pi(i)}\sse A^{(p)}_{\pi(i)}$, with the same
$c(S_i)=\lambda_2(A^{(i)})=k_0$.
The $p$ sets $S_i$ are $\s$-disjoint, since they are contained in
distinct parts of~$A^{(p)}$, which is itself $\s$-disjoint, 
but they all get the 
same color $\lambda_2(A^{(i)})$: contradiction.
\end{proof}

\begin{proof}[\bfseries Reduction of Theorem~\ref{t:main} to the case when
\boldmath$r$ is prime. ]
We proceed by induction on $r$, where we assume that the
result is true when $r$ is prime.
Thus let $\SS\sse 2^{[n]}$, let $r=r'r''$ with $2\le r',r''<r$, 
let $\chizero:=\chi(\KG^r\subds\SS)$, and assume that
\begin{equation}\tag{$*$}
\cd^r\subds\SS\ \ >\ \ (r-1)\chizero.
\end{equation}
We construct an auxiliary hypergraph $\TT\sse 2^{[n]}$
(on the same ground set as~$\SS$) by
\[
\TT\ \ :=\ \ \{N\sse[n]:\cd^{r'}\SS|_N>(r'-1)\chizero\}.
\]
Note that for this we use ``disjoint'' colorability defect,
corresponding to $\s=(1,...,1)$.
Using induction and the definition of~$\TT$, we now get
\[
(r'-1)\,\chi(\KG^{r'}\SS|_N)
\ \ \ge\ \ 
\cd^{r'}\SS|_N
\ \ >\ \ 
(r'-1)\chizero,
\]
and thus 
\begin{equation}\tag{$1$}\hspace{12mm}
\chi(\KG^{r'}\SS|_N)\ \ >\ \ \chizero\qquad\mbox{\rm for each $N\in\TT$.}
\end{equation}
\emph{Claim}: $\cd^{r''}\subds(\TT)>(r''-1)\chizero$.
\\[1mm]
\emph{Proof of the Claim}. Otherwise we could find an $\s$-disjoint
$r''$-family $N_1,\ldots, N_{r''}\sse[n]$ such that no $N_j$ contains
a set from~$\TT$ 
and such that $\sum_{j=1}^{r''}|N_j|\ge\ssum-(r''-1)K$.
In particular, none of the sets $N_j$ lies in~$\TT$, so by definition
of~$\TT$ we have
$\cd^{r'}(\SS|_{N_j})\le(r'-1)\chizero$ for all~$j$.
Thus for each~$j$ we can find $r'$ disjoint
sets $M_{j1},\ldots,M_{jr'}\sse N_j$, such that no~$M_{jk}$ contains
a set from~$\SS$, with $\sum_{k=1}^{r'}|M_{jk}|\ge|N_j|-(r'-1)\chizero$.

Taking all the sets $M_{jk}$ together, we have $r''r'=r$ 
subsets of~$[n]$, none of which contains a set from~$\SS$,
and they are $\s$-disjoint: they form an $\s$-disjoint union of
disjoint families. We compute
\begin{eqnarray*}
      \sum_{j=1}^{r''}\sum_{k=1}^{r'}|M_{jk}| 
&\ge& \sum_{j=1}^{r''} |N_j|   \ \ -\ \ r''(r'-1)\chizero \\
&\ge& \ssum\ -\ (r''-1)\chizero\ \ -\ \ r''(r'-1)\chizero \\
& = & \ssum\ -\ (r - 1)\chizero,
\end{eqnarray*}
which contradicts~($*$). Thus we have established the Claim. 
\\[2mm]
Using induction, together with the Claim, we get
\[
(r''-1)\,\chi(\KG^{r''}\subds(\TT))
\ \ \ge\ \ 
\cd^{r''}\subds(\TT)
\ \  > \ \ 
(r''-1)\chizero,
\]
and thus 
\begin{equation}\tag{$2$}
\chi(\KG^{r''}\subds(\TT))\ \ >\ \ \chizero.\hspace{25mm}
\end{equation}
Now consider a coloring $c:\SS\lra[\chizero]$ of $\KG^r\subds\SS$ 
by~$\chizero$ colors.
By ($1$), in every set $N\in\TT$ we find $r'$ disjoint sets from $\SS|_N$
which from $c$ get the same color $i\in[\chizero]$. 
Using this, we construct a new coloring $c':\TT\lra[\chizero]$
which assigns to $N\in\TT$ {one of} the (possibly several) colors~$i$
which $c$ assigns to $r'$ disjoint sets in~$\SS|_N$.
By ($2$), there are $r''$ sets $N_j\in\TT$, which are $\s$-disjoint,
and which from $c'$ get the same color $i_0=c'(N_j)$.
Thus we have $r''r'=r$ sets $M_{jk}\in\SS$ with $M_{jk}\sse N_j$,
also $\s$-disjoint, that get from~$c$ the same color~$i_0=c(M_{jk})$.
This contradicts the definition of~$\chizero$ and~$c$.
\end{proof}

\section{Chain Complexes and the \boldmath{$\ZZ_p$}-Tucker
  Lemma}\label{Sec:Zp-Tucker}

For convenience, the following is phrased in terms of chain complexes --
however, the argument is entirely combinatorial resp.\ 
easy to combinatorialize, since  
no homology, not even rank considerations, appear. 
For all the technology needed, Munkres 
\cite[esp.~\S\S 12-13]{Mun} is an excellent reference.

We start with a brief review
of chain complexes and chain homotopies, also intended to fix notation.
Let $\K$ be a finite abstract simplicial complex. 
The \emph{chain complex} $\CC(\K)$ of~$\K$ is
\[\CC(\K):\quad
\ldots \lra C_3
\stackrel{\partial_3}{\lra} C_2
\stackrel{\partial_2}{\lra} C_1
\stackrel{\partial_1}{\lra} C_0
\stackrel{\partial_0}{\lra}\{0\}\lra\ldots,
\]
where $C_k$ is the free abelian group of all formal linear
combinations of oriented $k$-faces of~$\K$, with integral
coefficients, and the boundary operators $\partial_k$ satisfy
$\partial_k\partial_{k+1}=0$. These are given by
$\partial_k[v_0,\ldots,v_k]=\sum^k_{i=0}(-1)^i[v_0,\ldots,\widehat{v_i}, 
\ldots,v_k]$.

A \emph{chain map} $\nu:\CC(\K)\lra \CC(\L)$
is a collection of homomorphisms $\nu_k: C_k(\K)\lra C_k (\L)$ 
such that $\partial \nu=\nu\partial$, that is, 
$\partial^{\sf L}_k \nu_k= \nu_{k-1} \partial^{\sf K}_k$
for all~$k$. Every simplicial map $f$ induces 
a chain map $f_\sharp=(f_{\sharp k})_{k\ge0}$.
Barycentric subdivision induces a canonical chain map~$\sd$.
Furthermore, any composition of chain maps is a chain map.

A \emph{chain homotopy} $D$ is a collection of homomorphisms
$D_k:C_k(\K)\lra C_{k+1}(\L)$, for all~$k$, 
with no compatibility condition. 
$\partial^{\sf L} D+D \partial^{\sf K}$ is then automatically a chain map. 
If  $\partial D+D\partial = \nu-\mu$, 
then $D$ is a chain homotopy \emph{between $\nu$ and $\mu$}.

If $\nu:\CC(\K)\lra\CC(\K)$ is a chain self-map, then its 
\emph{Lefschetz number} $\Lambda(\nu)$ is 
\[
\Lambda(\nu)\ \ :=\ \ \sum_k (-1)^k\,{\rm trace}(\nu_k).
\]
This {\bf counts} the nonempty
simplices that are mapped to themselves according
to the parity of their dimension and according to their effect on the
orientation.
For example, if $a:\K\lra \K$ is a constant map to
a vertex $v_0\in \K$, then $\Lambda(a_\sharp)=1$.
For ${\rm id}:\K\lra \K$, $\Lambda({\rm id}_\sharp)$ 
is the {Euler characteristic} of~$\K$.

\begin{proof}[\bfseries Combinatorial proof of the
\boldmath$\ZZ_p$-Tucker lemma. ]
Let us assume that, for some $m$, a $\ZZ_p$-equivariant coloring
$\lambda$ exists that does not produce a ``fully colored chain of
$p$ signed faces,'' as promised by the $\ZZ_p$-Tucker lemma.
\\
(1) In the first half of the proof we will
construct, under the assumption $\ssum\ge m(p-1)$, a
square of $\ZZ_p$-equivariant chain maps
\[
\begin{array}{ccc}
\CC(\sd\,(\sigma^{n-1})^{*p}\subds) &
\stackrel{\textstyle\lambda_\sharp}\longrightarrow & 
\CC(     (\sigma^{p-1}_{p-2})^{*m}) 
\\[2mm]
 {\sd}\Big\uparrow  &&
             \Big\downarrow \sd
\\[2mm]
\CC(     (\sigma^{n-1})^{*p}\subds) &
\stackrel{\textstyle \kappa_\sharp}\longleftarrow &
\CC(\sd\,(\sigma^{p-1}_{p-2})^{*m}).
\end{array}
\]
Here the vertices of $\sd(\sigma^{n-1})^{*p}\subds$
correspond to $0/1$-matrices of size $n\times p$ with row sums
at most $s_i$, as discussed above. The faces of the complex
correspond to chains of such matrices, with respect to the entrywise
$\le$-partial order. The $\ZZ_p$-action is free for prime~$p$, if $1\le s_i<p$.

The faces of $(\sigma^{p-1}_{p-2})^{*m}$ are
$0/1$-matrices of size $m\times p$ with no full row of ones:
we interpret them as admissible color sets. 
Again $\ZZ_p$ acts cyclically on the columns; this is free for 
prime~$p$. The map $\lambda$
of the $\ZZ_p$-Tucker lemma yields a simplicial map, and thus
the chain map $\lambda_\sharp$ used here.
The map is equivariant, by assumption.

The barycentric subdivision operators $\sd$, which yield
the vertical arrows in the square above, have 
explicit combinatorial descriptions that we do not
have to work out here. They are $\ZZ_p$-equivariant.

Finally, $\kappa: \sd(\sigma^{p-1}_{p-2})^{*m} \lra (\sigma^{n-1})^{*p}\subds$
is a simplicial map that we construct orbitwise, as follows.
The space to be mapped, $\sd(\sigma^{p-1}_{p-2})^{*m}$, 
is the barycentric subdivision of a simplicial 
complex of dimension $m(p-1)-1$, so it
is the order complex of a graded poset $Q_{m(p-1)}$ with $m(p-1)$ rank levels,
$\sd(\sigma^{p-1}_{p-2})^{*m}=\Delta(Q_{m(p-1)})$.
The free $\ZZ_p$-action on it respects the grading, so it decomposes
the rank levels of the poset $Q_{m(p-1)}$ into disjoint orbits of size~$p$.
The target space is a simplicial complex whose vertex set
is identified with the positions in an $n\times p$ matrix;
its faces are the $0/1$-matrices with at most $s_i$ ones in
the $i$-th row, for all $i$; the $\ZZ_p$-orbits of its vertices 
are exactly the rows of the matrix.
An equivariant simplicial map can now be defined orbitwise,
where the image of any element of an orbit determines
the images for all others.
We construct $\kappa$ such that the lowest $s_1$ rank levels are
mapped to  the lowest row of the $n\times p$ matrix. The next $s_2$ 
rank levels are mapped to the second lowest row of the matrix, etc.
Thus the orbits in the $\ell$-th rank level of $Q_{m(p-1)}$ are mapped
to the row number $\min\{t:\ell\le\sum_{i=1}^t s_i\}$.
(The following figure illustrates this for $n=4$, $p=5$, $s_i=2$.)
\[
\begin{picture}(0,0)%
\includegraphics{levelmap2.pstex}%
\end{picture}%
\setlength{\unitlength}{1973sp}%
\begingroup\makeatletter\ifx\SetFigFont\undefined%
\gdef\SetFigFont#1#2#3#4#5{%
  \reset@font\fontsize{#1}{#2pt}%
  \fontfamily{#3}\fontseries{#4}\fontshape{#5}%
  \selectfont}%
\fi\endgroup%
\begin{picture}(9923,2745)(601,-4659)
\put(10426,-2461){\makebox(0,0)[lb]{\smash{\SetFigFont{12}{14.4}{\rmdefault}{\mddefault}{\updefault}
$\Delta(Q_{m(p-1)})$
}}}
\put(601,-2386){\makebox(0,0)[lb]{\smash{\SetFigFont{12}{14.4}{\rmdefault}{\mddefault}{\updefault}
$(\sigma^{n-1})^{*p}\subds$
}}}
\end{picture}

\]
This is well-defined if the target matrix has enough rows,
that is, if $m(p-1)\le\sum_{i=1}^ns_i=\ssum$.
The definition on the vertices indeed 
yields a simplicial map into the target space:
any chain in $Q_{m(p-1)}$ contains at most
$s_i$ elements in the (at most $s_i$) adjacent 
rank levels that are mapped to the
$i$-th row of the $n\times p$ matrix.

In summary, for $\ssum\ge m(p-1)$, we can combine the four $\ZZ_p$-equivariant
chain maps of the square into a chain self-map
\[
\nu\ =\ \kappa_\sharp \sd\lambda_\sharp \sd:\quad
\CC((\sigma^{n-1})^{*p}\subds)\ \ \lra\ \ \CC((\sigma^{n-1})^{*p}\subds).
\]
Furthermore, the chain maps involved  are induced either by
simplicial maps, or by barycentric subdivision.
Thus all four of them, and thus in particular
$\nu$, are \emph{augmentation preserving} in the sense that they preserve
the sum of the coefficients of the vertices. 
\\
(2)
In the second half of the proof, we compute the Lefschetz number
of~$\nu$ in two ways. First, the $\ZZ_p$-actions are free
and the chain maps are $\ZZ_p$-equivariant, hence
the Lefschetz number of~$\nu$ satifies
\[
\Lambda(\nu)\equiv 0 \pmod p.
\]
However, we will show that if $\ssum\ge m(p-1)+1$, then $\nu$ 
restricts to an augmentation
preserving chain map of the chain complex of a cone, and thus necessarily has
\[
\Lambda(\nu)=1,
\]
which yields a contradiction for $m\le\frac{\ssum-1}{p-1}$.

If $\ssum\ge m(p-1)+1$, then we can extend $\kappa$ to a simplicial
map $\hat\kappa:\Delta(Q_{m(p-1)}\cup\{\hat1\})\lra(\sigma^{n-1})^{*p}\subds$,
where $\hat1$ denotes a new top element that is added to the poset
$Q_{m(p-1)}$. Indeed, just map this new element into the ``top row''
of the matrix; this gives a well-defined simplicial map
(not $\ZZ_p$-equivariant, of course).
The order complex $\Delta(Q_{m(p-1)}\cup\{\hat1\})$
is a cone, and thus so is its image 
$\K:=\hat\kappa(\Delta(Q_{m(p-1)}\cup\{\hat1\}))\sse(\sigma^{n-1})^{*p}\subds$:
the image of a cone under a simplicial map is always a cone.
We conclude that the image of~$\nu$ is contained in $\CC(\K)$,
where the restriction of~$\nu$ to $\CC(\K)$ has the
same Lefschetz number as $\nu$ itself. 
The following lemma thus completes the proof.
\end{proof}

(Barycentric subdivisions, as used in this
proof, do not yield simplicial maps; that's why the use of
chain complexes for this proof is essential.
The two barycentric subdivision operations that appear
in the square may be taken as a measure of 
complexity for the proof; in that sense, the proof of Schrijver's Theorem
given below is more complex; it needs $n-d+1$ barycentric subdivisions.)

\begin{lemma}\label{l:cone}
Let $\K=\K'* v_0$ be a finite simplicial cone, and let 
$\nu:\CC (\K)\longrightarrow \CC(\K)$ 
be an augmentation preserving chain map. Then $\Lambda(\nu)=1$.
\end{lemma}

\proof
The following five simple observations combine into a proof.
\\
(1) The identity map ${\rm id}:\K\longrightarrow \K$ and the
constant map to the apex 
$a:\K\longrightarrow\{v_0\}\sse \K$ 
are simplicial maps that induce chain maps
${\rm id}_\sharp, a_\sharp:\CC(\K)\longrightarrow\CC(\K)$.
These are chain homotopic: An explicit chain homotopy 
$D$, with $D_k:\CC_k(\K) \longrightarrow \CC_{k+1}(\K)$ for $k\ge0$,
is given by
\[
D:\sigma\longmapsto
\begin{cases}
v_0\ast\sigma & \text{if $v_0\notin\sigma$,}\\
0 & \text{otherwise.}
\end{cases}
\]
(2) If $\nu:\CC (\K)\longrightarrow \CC(\K)$ 
is {\em any} augmentation preserving chain map, then
\[
{\rm id}_\sharp \nu=\nu\qquad {\rm and}\qquad
a_\sharp \nu=a_\sharp.
\]
The first equality is clear, the second one is equivalent to being
augmentation preserving: $a_{\sharp0}$ maps every $0$-chain to ``sum of
coefficients times $[v_0]$,'' so we need that $\nu$ preserves
``sum of coefficients.'' 
\\
(3) $\bar D:= D \nu$ is a chain homotopy between
$\nu$ and $a_\sharp$. Indeed, using $\partial \nu=\nu\partial$ 
(since $\nu$ is a chain map), we get
\[
\nu-a_\sharp = 
{\rm id }_\sharp \nu-a_\sharp \nu
 = ({\rm id }_\sharp-a_\sharp) \nu
 = (\partial D + D\partial)\nu
 = \partial  (D \nu)+(D \nu) \partial.
\]
(4)
If two chain self-maps are connected by a chain homotopy, then they
have the same Lefschetz number. Indeed, let $\bar D$ be the chain homotopy,
then we compute
\begin{eqnarray*}
\Lambda(\partial\bar D+\bar D\partial) &=&
\sum_k (-1)^k\big[{\rm trace}(\partial_{k+1}\bar D_k)+
                  {\rm trace}(\bar D_{k-1}\partial_k)\big]\\ &=& 
\sum_k (-1)^k\big[{\rm trace}(\partial_{k+1}\bar D_k)+
                  {\rm trace}(\partial_k\bar D_{k-1})\big] ,
\end{eqnarray*}
which is a telescope sum that vanishes.
\\
(5) $\Lambda(a_\sharp)=1$.
\endproof

\section{A Special Case}\label{sec:badcolor}

Consider the case of the complete $k$-uniform
hypergraph $\SS={[n]\choose k}$, and of constant $\s=(s,...,s)$,
with $\ssum=ns$.
We get a lower bound on the chromatic number
of $\KG^r\subds{[n]\choose k}$ for this case by combination
of Theorem~\ref{t:main} with Lemma~\ref{l:cdefects},
and an upper bound from Lemma~\ref{l:col}:

\begin{cor}\label{cor:bounds}
Let $n\ge k\ge2$, $r>s\ge1$ with $r(k-1)\le ns$. Then
\[
1+\big\lceil\tfrac{ns-rk+1}{r-1}\big\rceil
\ \ \le\ \ \chi(\KG^r_{(s,...,s)}\textstyle{[n]\choose k})\ \ \le\ \ 
1+\big\lceil\tfrac1{\left\lfloor\tfrac{r-1}s\right\rfloor}\tfrac{ns-rk+1}s\big\rceil.
\]
In particular, if $s$ divides $r-1$, then
\[
\chi(\KG^r_{(s,...,s)}\textstyle{[n]\choose k})\ \ =\ \ 1+\big\lceil\tfrac{ns-rk+1}{r-1}\big\rceil.
\]
\end{cor}
This solves the generalized Kneser problem in the case when $\frac{r-1}s$
is an integer, which is the case, in particular, if $s=1$ (the
Alon-Frankl-Lov\'asz case).
However, if \mbox{$s\not|\, r-1$}, then --- contrary to the claim in
\cite[(3.3)]{Sa1} --- we do not have matching lower and upper bounds.
The smallest admissible parameters for this effect are
$s=2$, $r=4$, and $k=2$, and the following analysis shows that
in this case neither the lower bound nor the upper bound 
of Corollary \ref{cor:bounds} are sharp.

\begin{example}
Let $n\ge4$.
The generalized Kneser hypergraph $\KG^4_{(2,...,2)}{[n]\choose2}$
has the vertex set $E(K_n)$ (that is, the edges of a complete 
graph on $n$ vertices), while its edges are the
$4$-tuples of edges of~$K_n$
that form a subgraph of maximal degree at most~$2$.

The maximal sets $E'\sse E(K_n)$ that
don't contain an edge of $\KG^4_{(2,...,2)}{[n]\choose2}$
are of two types: either they have a vertex of degree at least~$3$, 
then they are ``a star plus one edge'' ($K_{1,n-1}+e$),
or every have no such vertex, then they
consist of exactly three disjoint edges.
Every $K$-coloring can be extended to 
a covering of $E(K_n)$ with $K$ such maximal sets.
From this one obtains that there are optimal colorings
for which the $i$-th color class is a $K_{1,n-i}+e$, and thus derives
that the chromatic number of $\KG^4_{(2,...,2)}{[n]\choose2}$ is
\[
\chi(\KG^4_{(2,...,2)}\textstyle{[n]\choose2})\ \ =\ \ 
\min\{K:\sum_{i=1}^K(n-i+1)\ge{n\choose2}\}\ \ =\ \ 
n-\big\lfloor\sqrt{2n+\tfrac14}-\tfrac12\big\rfloor.
\]
Thus, the chromatic number is roughly $n-\sqrt{2n}$ for
$r=4$ and $s=2$. This is to be compared with the lower bound
of Theorem~\ref{t:main}, which is 
$1+\lceil\frac{2n-7}3\rceil \approx \frac23n$, 
\emph{not very good}, 
and with the upper bound of Lemma~\ref{l:col}, which evaluates
to $1+\lceil\frac{2n-7}2\rceil=n-2$, \emph{useless}.
Thus, in the case where $r-1$ is not divisible by~$s$,
there is quite a gap between the upper and lower bounds 
in Corollary~\ref{cor:bounds}.
\end{example}

\section{Cyclic Oriented Matroids and Schrijver's Theorem}\label{sec:schrijver}

\begin{thm}[Schrijver \cite{Schrijver}]
For $n\ge2k>0$,
\[\textstyle
\chi(\KG^2{[n]\choose k}\stab)\ \ =\ \ n-2k+2.
\]
\end{thm}

This result is a strengthening of Lov\'asz' theorem: In the chain
\[\textstyle
n-2k+2\ \ \le\ \ 
\chi(\KG^2{[n]\choose k}\stab)\ \ \le\ \ 
\chi(\KG^2{[n]\choose k})\ \ \le\ \ n-2k+2,
\]
the first inequality is what we have to prove now,
the second one is trivial (coloring an induced subgraph), and
the third one is given by Kneser's coloring (Lemma~\ref{l:col}).
Schrijver \cite{Schrijver} indeed verified elegantly
that $\KG^2{[n]\choose k}\stab$
is a vertex-critical subgraph of the Kneser graph $\KG^2{[n]\choose k}$.
It is not edge-critical in general, as one may observe for
$k=2$, $n=3$, or less trivially for $k=2$, $n=6$. 

Let us also note that Schrijver's theorem is not implied by 
Do\lsoft nikov's, since Lemma~\ref{l:cdefects} provides smaller
(``worse'') values for $\cd^2{[n]\choose k}\stab$ than 
for~$\cd^2{[n]\choose k}$.

\begin{proof}[\bfseries Combinatorial proof of Schrijver's theorem. ]
Assume that we have a coloring
\[\textstyle
c:{[n]\choose k}\stab\ \ \lra\ \ [n-2k+1],
\]
and set $d:=n-2k+1$, so that $n>d\ge1$.
For $n\ge d$ we define 
\[
\Sigma^{d-1}(n)\ \ :=\ \ 
\Delta\big(\{X\in\pmz^n:\alt(X)\ge n-d+1\}\,,\,\le\big),
\]
the simplicial complex of all chains of sign vectors that
have an alternating subsequence with more than $n-d$ components.
It has a free $\ZZ_2$-action, given by $X\longleftrightarrow-X$.

The coloring $c$ of~${[n]\choose k}\stab$ now yields a simplicial,
$\ZZ_2$-equivariant map
\begin{eqnarray*}
\widehat{c}:\qquad
\Sigma^{d-1}(n)&\lra& \Sigma^{d-2}(d)\\
     (X^+,X^-) &\lmt& \textstyle
\big(c({X^+\choose k}\stab),c({X^-\choose k}\stab)\big).
\end{eqnarray*}
If $c$ is a correct coloring, then
the color sets $c({X^+\choose k}\stab)$ and $c({X^-\choose k}\stab)$ 
are disjoint. Moreover, $(X^+,X^-)\in\Sigma^{d-1}(n)$
has $\alt(X^+,X^-)\ge n-(d-1)=2k$, so
both $X^+\choose k$ and~$X^-\choose k$ contain at least one stable
$k$-set, so $c({X^+\choose k}\stab)$ and $c({X^-\choose k}\stab)$
cannot be empty, thus
$\alt(c({X^+\choose k}\stab),c({X^-\choose k}\stab))\ge2$, 
hence $\widehat c$ is well-defined for every vertex.
Furthermore, any chain of sign vectors is mapped by $\widehat c$ to a (weak)
chain of signed color sets, so we obtain a simplicial map.
This map is equivariant. 

\Topologically, $\Sigma^{d-1}(n)$ is a simplicial ($d-1$)-sphere,
namely the barycentric subdivision of the 
topological representation~\cite{FolkmanLawrence}
of the alternating oriented matroid 
of rank $d$ on $n$ elements, see \cite[Ex.~3.8]{BlandLasVergnas}
\cite[Chap.~5/Sect.~9.4]{BLSWZ} \cite{Z-hbo}.
Similarly, $\Sigma^{d-2}(d)$ is a simplicial ($d-2$)-sphere.
Both spheres have natural antipodal actions, and the map $\widehat c$
respects these. 
Thus, the Borsuk-Ulam theorem completes a topological proof at this
point, but we keep going on the combinatorial track.
\TopologicallyEnd

\noindent
Our next step is a quite trivial simplicial map,
\begin{eqnarray*}
\delta:\qquad
\Sigma^{d-2}(d)&\lra& \Sigma^{d-2}(d-1)\\
     (Y^+,Y^-) &\lmt& (Y^+\sm\{d\},Y^-\sm\{d\}),
\end{eqnarray*}
which deletes the last component of each sign vector $Y$.
This map is well-defined: deleting the last component
reduces $\alt(Y)$ at most by~$1$, the operation is compatible
with the partial order, and it is $\ZZ_2$-equivariant.

Now we use the canonical simplicial embedding maps of
$\Sigma^{d-2}(d-1)$ into a cone 
\[
\cone\Sigma^{d-2}(d-1)\ \ :=\ \ \Sigma^{d-2}(d-1)*\{v^+\}
\]
and then into a suspension 
\[
\susp\Sigma^{d-2}(d-1)\ \ :=\ \ \Sigma^{d-2}(d-1)*\{v^+,v^-\}.
\]
Thus we have maps
\[
\Sigma^{d-2}(d-1)
\ \ \stackrel{i'}{\lra}\ \ 
\cone\Sigma^{d-2}(d-1)
\ \ \stackrel{i''}{\lra}\ \ 
\susp\Sigma^{d-2}(d-1),
\]
where the composition $i:=i'' i'$ is $\ZZ_2$-equivariant
(with the natural $\ZZ_2$-action on the suspension that
interchanges $v^+$ and~$v^-$).

Finally, we will construct an equivariant chain map
\[
\xi:\ \CC(\susp\Sigma^{d-2}(d-1))\ \ \lra\ \ \CC(\Sigma^{d-1}(d)),
\]
and sequence of equivariant chain maps
\[
\zeta^\ell:\ \CC(\Sigma^{d-1}(\ell-1))\ \ \lra\ \ \CC(\Sigma^{d-1}(\ell)),
\]
for $d<\ell\le n$. Once these maps have been 
constructed, the proof will be complete, since then we
have a square of equivariant, augmentation preserving
chain maps
\[
\begin{array}{rcc}
\CC(\Sigma^{d-1}(n)) &
\stackrel{\textstyle \widehat c_\sharp }\longrightarrow & 
\CC(\Sigma^{d-2}(d)) 
\\[2mm]
 {\zeta^n}\Big\uparrow \qquad && \Big\downarrow \delta_\sharp
\\[2mm]
\vdots\qquad\ && \CC(\Sigma^{d-2}(d-1))  
\\[2mm]
 {\zeta^{d+1}}\Big\uparrow \qquad  &&
             \Big\downarrow i_\sharp
\\[2mm]
\CC(\Sigma^{d-1}(d)) &
\stackrel{\textstyle \xi}\longleftarrow & 
\CC(\susp\Sigma^{d-2}(d-1)) .
\end{array}
\]
Indeed, as in Section~\ref{Sec:Zp-Tucker}
we can then argue that the composition
\[
i_\sharp\delta_\sharp\widehat c_\sharp\zeta^n\zeta^{n-1}\cdots\zeta^{d+1}\xi:\ 
\CC(\susp\Sigma^{d-2}(d-1))\ \ \lra\ \ \CC(\susp\Sigma^{d-2}(d-1))
\]
is $\ZZ_2$-equivariant, so it has even Lefschetz number,
but it also restricts to a cone (the image of $i'_\sharp$ is contained
in $\CC(\cone\Sigma^{d-2}(d-1))$), and thus its Lefschetz number
is~$1$ (Lemma~\ref{l:cone}).

The chain maps $\xi$ and $\zeta^\ell$ (for $d+1\le\ell\le n$)
can be written down combinatorially, 
by giving a formula for the image of an arbitrary
$k$-simplex as a sum of $k$-simplices with $\pm1$-coefficients. 
However, for the exposition we prefer to give a geometric
description, from which the combinatorial one can then be derived.

For $\xi$, note that 
$\Sigma^{d-1}(d)$ can be interpreted as the boundary complex
of~$\sd([-1,1]^d)$, whose ``equator'' subsphere naturally corresponds to 
$\Sigma^{d-2}(d-1)\cong \sd([-1,1]^{d-1})$.
\[
\includegraphics[height=4.6cm]{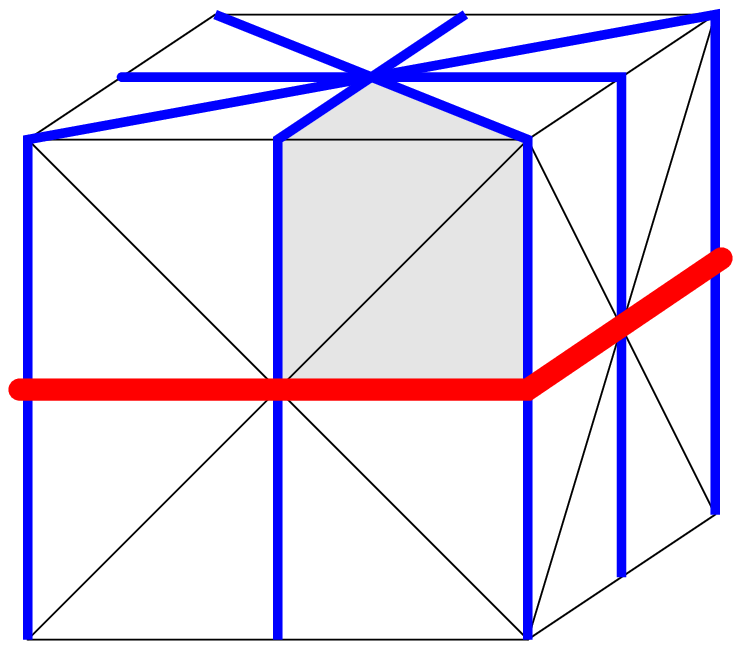}
\]
This suggests a natural subdivision chain map
$\CC(\susp\Sigma^{d-2}(d-1))\lra\CC(\Sigma^{d-1}(d))$, 
as indicated in the figure, where every
simplex in $\Sigma^{d-2}(d-1)$ is mapped ``to itself'' on the equator, 
the vertices $v^\pm$ are mapped to the north/south poles $\pm e_d$, 
and any $k$-simplex $\sigma*v^+$, say, is mapped to a signed sum of
$k+1$ $k$-simplices, $k$ of which triangulate the prism
$\sigma\times[0,1]$, while the last one is $+\,\sigma{\times}\{1\}*e_d$.
\medskip

The construction of the chain maps 
$\zeta^\ell:\CC(\Sigma^{d-1}(\ell-1))\lra\CC(\Sigma^{d-1}(\ell))$
is quite similar:
$\Sigma^{d-1}(\ell-1)$ is the face poset of the
dual cyclic oriented matroid of rank $d$ on $\ell-1$ elements,
as described and analyzed in~\cite{Z-hbo}.

The following figure illustrates the sign vectors associated to various
faces (vertices and edges) of $\Sigma^{d-1}(\ell-1)$, for
$\ell=5$ and $d=3$.
\[
\begin{picture}(0,0)%
\includegraphics{cyclic1.pstex}%
\end{picture}%
\setlength{\unitlength}{1973sp}%
\begingroup\makeatletter\ifx\SetFigFont\undefined%
\gdef\SetFigFont#1#2#3#4#5{%
  \reset@font\fontsize{#1}{#2pt}%
  \fontfamily{#3}\fontseries{#4}\fontshape{#5}%
  \selectfont}%
\fi\endgroup%
\begin{picture}(6032,6495)(585,-5986)
\put(2701,-3286){\makebox(0,0)[lb]{\smash{\SetFigFont{12}{14.4}{\rmdefault}{\mddefault}{\updefault}
$3$
}}}
\put(5326,-1036){\makebox(0,0)[lb]{\smash{\SetFigFont{12}{14.4}{\rmdefault}{\mddefault}{\updefault}
$1$
}}}
\put(2251,-586){\makebox(0,0)[b]{\smash{\SetFigFont{12}{14.4}{\rmdefault}{\mddefault}{\updefault}
$2$
}}}
\put(4126,-4636){\makebox(0,0)[b]{\smash{\SetFigFont{12}{14.4}{\rmdefault}{\mddefault}{\updefault}
$4$
}}}
\put(1951,-2011){\makebox(0,0)[lb]{\smash{\SetFigFont{8}{9.6}{\rmdefault}{\mddefault}{\updefault}
${+}00{-}$
}}}
\put(3376,-4186){\makebox(0,0)[lb]{\smash{\SetFigFont{8}{9.6}{\rmdefault}{\mddefault}{\updefault}
${+}{-}00$
}}}
\put(5326,-5536){\makebox(0,0)[lb]{\smash{\SetFigFont{8}{9.6}{\rmdefault}{\mddefault}{\updefault}
$0{-}0{+}$
}}}
\put(2626,-4861){\makebox(0,0)[lb]{\smash{\SetFigFont{8}{9.6}{\rmdefault}{\mddefault}{\updefault}
${+}0{-}{+}$
}}}
\put(4426,-4936){\makebox(0,0)[lb]{\smash{\SetFigFont{8}{9.6}{\rmdefault}{\mddefault}{\updefault}
${+}{-}0{+}$
}}}
\put(1951,-3961){\makebox(0,0)[rb]{\smash{\SetFigFont{8}{9.6}{\rmdefault}{\mddefault}{\updefault}
${+}0{-}0$
}}}
\put(2251,-1186){\makebox(0,0)[lb]{\smash{\SetFigFont{8}{9.6}{\rmdefault}{\mddefault}{\updefault}
${+}0{+}{-}$
}}}
\put(3751,-2311){\makebox(0,0)[lb]{\smash{\SetFigFont{8}{9.6}{\rmdefault}{\mddefault}{\updefault}
${+}{-}{+}{-}$
}}}
\put(3001,314){\makebox(0,0)[lb]{\smash{\SetFigFont{8}{9.6}{\rmdefault}{\mddefault}{\updefault}
$00{+}{-}$
}}}
\put(3451,-5986){\makebox(0,0)[lb]{\smash{\SetFigFont{8}{9.6}{\rmdefault}{\mddefault}{\updefault}
$00{-}{+}$
}}}
\end{picture}

\]
The vertex set of $\Sigma^{d-1}(\ell-1)$ is
$\{Y\in\pmz^\ell:\alt(Y)\ge \ell-d\}$. 
Thus $\Sigma^{d-1}(\ell-1)$ decomposes
into a \emph{positive hemisphere} $\Gamma^+$, the induced subcomplex
given by all vertices with $\alt(Y,+)>\ell-d$, and 
the \emph{negative hemisphere} $\Gamma^-$
of all vertices with $\alt(Y,-)>\ell-d$. Every simplex is contained
in one of these two ``hemispheres,'' whose intersection is
the \emph{equator}~$\Gamma^0$, induced on all vertices $Y$ with
$\alt(Y)>\ell-d$. The equator is naturally isomorphic to
$\Sigma^{d-2}(\ell-1)$.

Similarly, the simplicial complex $\Sigma^{d-1}(\ell)$
decomposes into a positive hemisphere $\widehat\Gamma^+$,
on the vertices $X$ with last component $X_n\in\{{+},0\}$,
and a negative hemisphere $\widehat\Gamma^-$, 
on the vertices $X$ with last component $X_n\in\{{-},0\}$,
whose intersection is the equator $\widehat\Gamma^0$, 
given by the vertices $X$ with $X_n=0$, which is again isomorphic to  
$\Sigma^{d-2}(\ell-1)$.

For the construction of~$\zeta^\ell$,
let $\sigma$ be an (oriented) simplex of $\Sigma^{d-1}(\ell)$;
since we are dealing with the order complex of a poset, all
simplices have a natural ordering of their vertices, and thus
they have natural orientations.
We distinguish three cases. 
\\[3mm]
(1) If $\sigma$ lies in the equator $\Gamma^0$, then
$\zeta^\ell(\sigma)$ is the corresponding oriented simplex in
$\widehat\Gamma^0$.
\\[3mm]
(2) If $\sigma$ has more than one vertex that does not lie in the
equator, then it is mapped to the corresponding simplex in the
positive or negative hemisphere.
\\[3mm]
(3) Assume that $\sigma=[Y_0<Y_1<\ldots<Y_k]$ has a facet in the equator,
but does not itself lie in the equator; that is, all its vertices
lie in the equator, except for $Y_0$, which lies in the interior of
the positive hemisphere, say. Then $\sigma$ is mapped to
a sum of~$k+1$ oriented $k$-simplices in $\widehat\Gamma^+$,
which again correspond to 
$k$ simplices that triangulate a prism over~$\sigma$, plus one extra
simplex in the interior of the positive hemisphere, and which in the 
following formal description comes first:

\[
\begin{array}{rlr}
[Y_0<Y_1<\ldots<Y_k]&\ \longmapsto\\[2mm]
     &[(Y_0,{+})<(Y_1,{+})<(Y_2,{+})<&\ldots\quad <(Y_{k-1},{+})< (Y_k,{+})]\\
-    &[(Y_1, 0 )<(Y_1,{+})<(Y_2,{+})<&\ldots\quad <(Y_{k-1},{+})< (Y_k,{+})]\\
+    &[(Y_1, 0 )<(Y_2, 0 )<(Y_2,{+})<&\ldots\quad <(Y_{k-1},{+})< (Y_k,{+})]\\
            &\quad \vdots \qquad\qquad\vdots     &  \vdots         \qquad  \\
+(-1)^{k-1} &[(Y_1, 0 )<(Y_2, 0 )<\quad\ldots    &
                               \!\!  < (Y_{k-1},0)<(Y_{k-1},{+})<(Y_k,{+})]\\
+(-1)^k     &[(Y_1, 0 )<(Y_2, 0 )<\quad\ldots    &
                               \!\!  < (Y_{k-1},0)<(Y_{k-1},\,0\,)<(Y_k,{+})]
\end{array}
\]
Our figure tries to illustrate this for $d=3$ and $\ell=5$.
In the left figure,
which represents $\Sigma^{d-1}(\ell-1)=\Sigma^2(4)$, 
the equator arises as a subcomplex
of the barycentric subdivision. The right figure, depicting
$\Sigma^{d-1}(\ell)=\Sigma^2(5)$, has the equator as a regular subsphere.
$\zeta^5$ is a chain map from the left simplicial complex
to the right one; it maps the shaded triangle 
$[({+}0{-}0) < ({+}0{-}{+}) < ({+}{+}{-}{+})]$
on the left to the sum of the three triangles shaded on the right:
\\[2mm]
\begin{picture}(0,0)%
\includegraphics{cyclic2.pstex}%
\end{picture}%
\setlength{\unitlength}{2329sp}%
\begingroup\makeatletter\ifx\SetFigFont\undefined%
\gdef\SetFigFont#1#2#3#4#5{%
  \reset@font\fontsize{#1}{#2pt}%
  \fontfamily{#3}\fontseries{#4}\fontshape{#5}%
  \selectfont}%
\fi\endgroup%
\begin{picture}(12632,6032)(585,-5777)
\put(916,-4951){\makebox(0,0)[lb]{\smash{\SetFigFont{9}{10.8}{\rmdefault}{\mddefault}{\updefault}
${+}{+}{-}{+}$
}}}
\put(2791,-5236){\makebox(0,0)[rb]{\smash{\SetFigFont{9}{10.8}{\rmdefault}{\mddefault}{\updefault}
${+}0{-}{+}$
}}}
\put(1921,-3871){\makebox(0,0)[rb]{\smash{\SetFigFont{9}{10.8}{\rmdefault}{\mddefault}{\updefault}
${+}0{-}0$
}}}
\put(2701,-3286){\makebox(0,0)[lb]{\smash{\SetFigFont{12}{14.4}{\rmdefault}{\mddefault}{\updefault}
$3$
}}}
\put(5326,-1036){\makebox(0,0)[lb]{\smash{\SetFigFont{12}{14.4}{\rmdefault}{\mddefault}{\updefault}
$1$
}}}
\put(2251,-586){\makebox(0,0)[b]{\smash{\SetFigFont{12}{14.4}{\rmdefault}{\mddefault}{\updefault}
$2$
}}}
\put(4126,-4636){\makebox(0,0)[b]{\smash{\SetFigFont{12}{14.4}{\rmdefault}{\mddefault}{\updefault}
$4$
}}}
\end{picture}

\smallskip

Now one verifies either geometrically (depending on the geometric
realization of the spheres in question as barycentric subdivisions
of arrangement spheres), or combinatorially, that these rules satisfy 
$\partial\zeta^\ell\sigma=\zeta^{\ell-1}\partial\sigma$
for all $\sigma$, that is, they
provide equivariant simplicial chain maps $\zeta^\ell$, as required.
\end{proof}

\section{Stable Kneser Hypergraphs}\label{s:conj:big}

The stable Kneser hypergraphs $\KG^r{[n]\choose k}\tstab$ 
are induced subgraphs 
of the usual Kneser hypergraphs $\KG^r{[n]\choose k}$. 
For $t=r\ge2$, $k\ge2$ and $n\ge kr$, it seems that these
sub-hypergraphs have the same chromatic numbers as the full hypergraphs, 
\[\textstyle
\chi\big(\KG^r{[n]\choose k}\rstab\big)\ \ =\ \ 
\lceil\tfrac{n-(k-1)r}{r-1}\rceil.
\]
Here ``$\le$'' holds by Lemma~\ref{l:col}. Furthermore, 
if $r-1$ divides $n-rk$, 
then it appears that~$\KG^r{[n]\choose k}\rstab$
is vertex-critical, that is, for every $S\in{[n]\choose k}\stab$,
\[\textstyle
\chi\big(\KG^r{[n]\choose k}\stab\sm\{S\}\big)\ \ <\ \ 
\lceil\tfrac{n-(k-1)r}{r-1}\rceil.
\]
On the other hand, the $r$-stable $r$-th Kneser hypergraph is
not vertex critical in general, for example for $n-rk=1$ and~$r>2$.

\small
\subsection*{Acknowledgements}
Thanks to Jirka Matou\v sek for inspiring preprints, lectures, and discussions.
\\
Thanks to James Munkres for teaching me Algebraic Topology
the combinatorial way \cite{Mun}, and to
Anders Bj\"orner for teaching me Combinatorics the 
topological way \cite{Bj}.
\\
Thanks to Torsten Heldmann for special support in special times.

\bibliographystyle{alpha}
\bibliography{combi}

\end{document}